\documentclass[reqno,10pt]{amsart}

\usepackage{amsmath,amsfonts,amsthm,amssymb,amscd}
\usepackage[dvips]{graphics}
\usepackage{epsfig}

\newcommand{\ncr}[2]{{#1 \choose #2}}
\newcommand\be{\begin{equation}}
\newcommand\ee{\end{equation}}
\newcommand\bea{\begin{eqnarray}}
\newcommand\eea{\end{eqnarray}}
\newcommand\bi{\begin{itemize}}
\newcommand\ei{\end{itemize}}
\newcommand\ben{\begin{enumerate}}
\newcommand\een{\end{enumerate}}
\newcommand\bc{\begin{center}}
\newcommand\ec{\end{center}}
\newcommand\ba{\begin{array}}
\newcommand\ea{\end{array}}

\newtheorem{thm}{Theorem}[section]

\newtheorem{lem}[thm]{Lemma}

\newtheorem{defi}[thm]{Definition}

\newtheorem{rek}[thm]{Remark}

\newcommand{\sumin}{\sum_{i=1}^N}
\newcommand{\fon}{\frac{1}{N}}

\newcommand{\gl}{\lambda}

\newcommand{\R}{\ensuremath{\mathbb{R}}}

\newcommand{\N}{\mathbb{N}}

\newcommand{\E}{\ensuremath{\mathbb{E}}}
\newcommand{\p}{\ensuremath{\mathbb{P}}}
\newcommand{\f}{\mathcal{F}}

\newcommand{\foh}{\frac{1}{2}}  

\newcommand{\mattwo}[4]
{\left(\begin{array}{cc}
                        #1  & #2   \\
                        #3 &  #4
                          \end{array}\right) }

\newcommand{\ga}{\alpha}                  
\newcommand{\gb}{\beta}
\newcommand{\gep}{\epsilon}

\numberwithin{equation}{section}

\begin{document}

\title[Real
Symmetric Palindromic Toeplitz Matrices and Circulant Matrices\ \
]{Distribution of Eigenvalues of Real Symmetric Palindromic Toeplitz
Matrices and Circulant Matrices}

\author{Adam Massey}
\email{Adam\underline{\ }Massey@brown.edu}
\address{Department of Mathematics, Brown University, Providence,
RI 02912\bigskip\bigskip}

\author{Steven J. Miller}
\email{sjmiller@math.brown.edu}
\address{Department of Mathematics, Brown University, Providence,
RI 02912\bigskip\bigskip}

\author{John Sinsheimer}
\email{sinsheimer.2@osu.edu}
\address{Department of Mathematics, The Ohio State University,
Columbus, OH 43210}

\date{\today}



\begin{abstract}
Consider the ensemble of real symmetric Toeplitz matrices, each
independent entry an i.i.d. random variable chosen from a fixed
probability distribution $p$ of mean 0, variance 1, and finite
higher moments. Previous investigations showed that the limiting
spectral measure (the density of normalized eigenvalues) converges
weakly and almost surely, independent of $p$, to a distribution
which is almost the standard Gaussian. The deviations from Gaussian
behavior can be interpreted as arising from obstructions to
solutions of Diophantine equations. We show that these obstructions
vanish if instead one considers real symmetric palindromic Toeplitz
matrices, matrices where the first row is a palindrome. A similar
result was previously proved for a related circulant ensemble
through an analysis of the explicit formulas for eigenvalues. By
Cauchy's interlacing property and the rank inequality, this ensemble
has the same limiting spectral distribution as the palindromic
Toeplitz matrices; a consequence of combining the two approaches is
a version of the almost sure Central Limit Theorem. Thus our
analysis of these Diophantine equations provides an alternate
technique for proving limiting spectral measures for certain
ensembles of circulant matrices.
\end{abstract}

\maketitle

\tableofcontents


\section{Introduction}

\subsection{History}
Random matrix theory has successfully modeled many complicated
systems, ranging from energy levels of heavy nuclei in physics to
zeros of $L$-functions in number theory. For example, while the
nuclear structure of hydrogen is quite simple and amenable to
description, the complicated interactions of the over 200 protons
and neutrons in a Uranium nucleus prevent us from solving the
Hamiltonian equation (let alone even writing down the entries of the
matrix!). Similar to statistical mechanics, the complexity of the
system actually helps us describe the general features of the
solutions. Wigner's great insight was to approximate the infinite
dimensional Hamiltonian matrix with the limit of $N\times N$ real
symmetric matrices chosen randomly (each independent entry is chosen
from a Gaussian density; this ensemble of matrices is called the GOE
ensemble). For each $N$ one can calculate averages over the weighted
set of matrices, such as the density of or spacings between
normalized eigenvalues. Similar to the Central Limit Theorem, as
$N\to\infty$ with probability one we have that the behavior of the
normalized eigenvalues of a generic, randomly chosen matrix agrees
with the limits of the system averages.

Instead of choosing the entries of our matrices from Gaussian
densities, we could instead choose a nice probability distribution
$p$, for example, a distribution with mean $0$, variance $1$ and
finite higher moments. For real symmetric matrices with independent
entries i.i.d.r.v. from suitably restricted probability
distributions, the limiting distribution of the density of
normalized eigenvalues is the semi-circle density (see
\cite{Wig,Meh}). While there is universality in behavior of the
density of normalized eigenvalues, much less can be proved for the
distribution of normalized spacings; though extensive numerical
investigations support the conjecture that the behavior is the same
as the GOE ensemble, this cannot be proved for general $p$.

It is a fascinating question to impose additional structure on the
real symmetric matrices, and see how the behavior changes. The GOE
ensemble has $N(N+1)/2$ independent parameters, the $a_{ij}$ with $i
\le j \in \{1,\dots, N\}$. For sub-ensembles, different limiting
distributions arise. For example, to any graph $G$ one can associate
its adjacency matrix $A_G$, where $a_{ij}$ is the number of edges
connecting vertices $i$ and $j$. If $G$ is a simple $d$-regular
graph with no self-loops (there is at most one edge between two
vertices, each vertex is connected to exactly $d$ vertices, and
there are no edges from a vertex to itself), its adjacency matrix is
all $0$'s and $1$'s. Such graphs often arise in network theory. The
eigenvalues of these adjacency matrices are related to important
properties of the graphs: all eigenvalues lie in $[-d, d]$, $d$ is a
simple eigenvalue if and only if the graph is connected, and if the
graph is connected then the size of the second largest eigenvalue is
related to how quickly information propagates in the network (see,
for example, \cite{DSV}). Instead of choosing the matrix elements
randomly, for each $N$ there are only finitely many $N\times N$
$d$-regular graphs, and we choose uniformly from this set. While
$d$-regular graphs are a subset of real symmetric matrices, they
have different behavior. McKay \cite{McK} proved the density of
eigenvalues of $d$-regular graphs is given by Kesten's Measure, not
the semi-circle; however, as $d\to\infty$ the distributions converge
to the semi-circle density. Interestingly, numerical simulations
support the conjecture that the spacings between normalized
eigenvalues are the same as the GOE; see for example \cite{JMRR}.

Thus by examining sub-ensembles, one has the exciting possibility of
seeing new, universal distributions and behavior; for adjacency
matrices of $d$-regular graphs, only $dN/2$ of the possible
$N(N-1)/2$ edges are chosen, and the corresponding $a_{ij}$ (which
equal $1$) are the only non-zero entries of the adjacency matrices.
Recently the density of eigenvalues of another thin subset of real
symmetric matrices was studied. Recall an $N\times N$ Toeplitz
matrix $A_N$ is of the form \be A_N \ = \ \left(\begin{array}{ccccc}
b_{0}  &  b_{1}  & b_{2}  & \cdots & b_{N-1} \\
b_{-1} &  b_{0}  & b_{1}  & \cdots & b_{N-2} \\
b_{-2} &  b_{-1} & b_{0}  & \cdots & b_{N-3} \\
\vdots & \vdots  & \vdots & \ddots & \vdots \\
            b_{1-N} &  b_{2-N} & b_{3-N} & \cdots & b_0 \\
  \end{array}\right), \ \ \ \ a_{ij} \ = \
  b_{j-i}.
\ee Bai \cite{Bai} proposed studying the density of eigenvalues of
real symmetric Toeplitz matrices with independent entries
independently drawn from a nice distribution $p$. As a Toeplitz
matrix has $N$ degrees of freedom (the $b_i$'s), this is a very thin
sub-ensemble of all real symmetric matrices, and the imposed
structure leads to new behavior.

Initial numerical simulations suggested that the density of
normalized eigenvalues might converge to the standard Gaussian
density; however, Bose-Chatterjee-Gangopadhyay \cite{BCG} showed
this is not the case by calculating the fourth moment of the
limiting spectral measure (see Definitions \ref{defi:nesd} and
\ref{defi:lsd}) of the normalized eigenvalues. The fourth moment is
$2 \frac23$, close to but not equal to the standard Gaussian
density's fourth moment of $3$. Bryc-Dembo-Jiang \cite{BDJ}
(calculating the moments using uniform variables and interpreting
the results as volumes of solids related to Eulerian numbers) and
Hammond-Miller \cite{HM} (calculating the moments by solving systems
of Diophantine equations with obstructions) then independently found
somewhat intractable formulas for all the moments, and further
quantified the non-Gaussian behavior. The analysis in \cite{HM}
shows that the moments of the Toeplitz ensemble grow fast enough to
give a distribution with unbounded support, but significantly slower
than the standard Gaussian's moments (the ratio of the
$2k$\textsuperscript{th} Toeplitz moment to the standard Gaussian's
moment tends to zero as $k\to\infty$).

In \cite{HM} it was observed that their techniques may be applicable
to a related ensemble. Specifically, by imposing an additional
symmetry on the matrices by requiring that the first row be a
palindrome (see \eqref{eq:defrsptmat}), the obstructions to the
Diophantine equations vanish and the limiting spectral measure
converges weakly, in probability and almost surely to the standard
Gaussian (see \S\ref{sec:intromainresults} for the exact
statements). Bose and Mitra \cite{BM} proved weak convergence for an
ensemble closely related to our palindromic Toeplitz matrices (see
\eqref{eq:defsymmtoep} for the ensemble they studied). They combined
explicit expressions for the eigenvalues of circulant matrices and
probabilistic arguments to construct the empirical spectral
distribution; with these in place, they then show the limiting
spectral distribution is the standard Gaussian.

We show in Theorem \ref{thm:symmtoepispaltoep} that our analysis
gives an alternate proof of Bose and Mitra's result. We generalize
the linear algebra arguments described in \cite{HM} to analyze the
Diophantine equations that arise. The eigenvalues of our palindromic
Toeplitz ensemble are interlaced with those of the circulant
ensemble of \eqref{eq:defsymmtoep}. By Cauchy's interlacing property
(Lemma \ref{lem:cauchyinterlace}) and the rank inequality (Lemma
\ref{lem:rankineq}), our analysis of the Diophantine equations
related to the palindromic Toeplitz ensemble provides an alternate
proof of the limiting spectral measure of the circulant ensemble in
\eqref{eq:defsymmtoep}. This equivalence may be of use to other
researchers studying related ensembles, as we have replaced having
to calculate and work with explicit formulas for eigenvalues to
solving a system of Diophantine equations without obstructions.
Additionally, this equivalence leads to a version of the almost sure
Central Limit Theorem (see Theorem \ref{thm:asclt2}).

\subsection{Notation}\label{sec:subnotation}

We briefly review the notions of convergence examined in this paper
(see \cite{GS} for more details) and define the quantities studied.
We consider real symmetric palindromic Toeplitz matrices whose
independent entries are i.i.d. random variables chosen from some
distribution $p$ with mean 0, variance 1, and finite higher moments.
\emph{For convenience we always assume $N$ is even.} Thus our
matrices are of the form \be\label{eq:defrsptmat} A_N\ =\ \left(
\begin{array}{ccccccc}
b_0    &  b_1   &  b_2   &  \cdots  &  b_2   &  b_1  &  b_0  \\
b_1    &  b_0   &  b_1   &  \cdots  &  b_3   &  b_2  &  b_1  \\
b_2    &  b_1   &  b_0   &  \cdots  &  b_4   &  b_3  &  b_2  \\
\vdots & \vdots & \vdots &  \ddots  & \vdots & \vdots &  \vdots  \\
b_2    & b_3    & b_4    &  \cdots  & b_0    & b_1   &  b_2   \\
b_1    & b_2    & b_3    &  \cdots  & b_1    & b_0   &  b_1  \\
b_0    & b_1    & b_2    &  \cdots  & b_2    & b_1   &  b_0
  \end{array}\right).
\ee 
Each $N \times N$ matrix $A_N$ is parametrized by $N/2$ numbers:
$b_0(A_N),\dots,b_{N/2-1}(A_N)$. We may thus identify such $N\times
N$ real symmetric palindromic Toeplitz matrices with vectors in
$\R^{N/2}$.

For each integer $N$ let $\Omega_N$ denote the set of $N \times N$
real symmetric palindromic Toeplitz matrices. We construct a
probability space $(\Omega_N,\f_N,\p_N)$ by setting \bea & &
\p_N\left(\left\{A_N\in\Omega_N: b_{iN}(A_N) \in [\ga_i, \gb_i]\
{\rm for}\ i \in \left\{0,\dots,N/2-1\right\}\right\}\right)
\nonumber\\ & & \ \ \ \ \ \ \ \ = \ \prod_{i=1}^M
\int_{x_i=\ga_i}^{\gb_i} p(x_i)dx_i, \eea where each $dx_i$ is
Lebesgue measure. To each $A_N\in \Omega_N$ we attach a spacing
measure by placing a point mass of size $1/N$ at each normalized
eigenvalue\footnote{From the eigenvalue trace lemma
($\text{Trace}(A_N^2) = \sum_i \gl_i^2(A_N)$) and the Central Limit
Theorem, we see that the eigenvalues of $A_N$ are of order
$\sqrt{N}$. This is because $\text{Trace}(A_N^2) = \sum_{i,j=1}^N
a_{ij}^2$, and since each $a_{ij}$ is drawn from a mean $0$,
variance $1$ distribution, $\text{Trace}(A_N^2)$ is of size $N^2$.
This suggests the appropriate scale for normalizing the eigenvalues
is to divide each by $\sqrt{N}$.} $\gl_i(A_N)$: \be\label{eq:MkANl}
\mu_{A_N}(x)dx \ = \ \frac{1}{N} \sum_{i=1}^N \delta\left( x -
\frac{\gl_i(A_N)}{\sqrt{N}} \right)dx, \ee where $\delta(x)$ is the
standard Dirac delta function. We call $\mu_{A_N}$ the normalized
spectral measure associated to $A_N$.

\begin{defi}[Normalized empirical spectral
distribution]\label{defi:nesd} Let $A_N$ be an $N\times N$ real
symmetric matrix with eigenvalues $\lambda_N \ge \cdots \ge
\lambda_1$. The normalized empirical spectral distribution (the
empirical distribution of normalized eigenvalues) $F^{A_N/\sqrt{N}}$
is defined by \be F^{A_N/\sqrt{N}}(x) \ = \frac{\#\{i \le N:
\lambda_i/\sqrt{N} \le x\}}{N}. \ee \end{defi}

As $F^{A_N/\sqrt{N}}(x) = \int_{-\infty}^x \mu_{A_N}(t)dt$, we see
that $F^{A_N/\sqrt{N}}$ is the cumulative distribution function
associated to the measure $\mu_{A_N}$. 

We are interested in the behavior of a typical $F^{A_N/\sqrt{N}}$ as
$N\to\infty$. Our main results are that $F^{A_N/\sqrt{N}}$ converges
to the cumulative distribution function of the Gaussian (we describe
the type of convergence in \S\ref{sec:intromainresults}). Thus let
$M_m$ equal the $m$\textsuperscript{th} moment of the standard
Gaussian (so $M_{2k} = (2k-1)!!$ and $M_{2k+1} = 0$). As there is a
one-to-one correspondence between $N\times N$ real symmetric
palindromic Toeplitz matrices and $\R^{N/2}$, we may study the more
convenient infinite sequences. Thus our outcome space is $\Omega_\N
= \{b_0,b_1,\dots\}$, and if $\omega = (\omega_0,\omega_1,\dots) \in
\Omega_\N$ then \be {\rm Prob}(\omega_i \in [\alpha_i,\beta_i])\ =\
\int_{\alpha_i}^{\beta_i} p(x_i)dx_i.\ee We denote elements of
$\Omega_\N$ by $A$ to emphasize the correspondence with matrices,
and we set $A_N$ to be the real symmetric palindromic Toeplitz
matrix obtained by truncating $A = (b_0,b_1,\dots)$ to
$(b_0,\dots,b_{N/2-1})$. We denote the probability space by
$(\Omega_\N,\f_\N,\p_\N)$.

To each integer $m \ge 0$ we define the random variable $X_{m;N}$ on
$\Omega_\N$ by \be X_{m;N}(A) \ = \ \int_{-\infty}^\infty x^m
dF^{A_N/\sqrt{N}}(x); \ee note this is the $m$\textsuperscript{th}
moment of the measure $\mu_{A_N}$.

We investigate several types of convergence.

\ben

\item (Almost sure convergence) For each $m$, $X_{m;N} \to X_m$
almost surely if \be \p_\N\left(\{A \in \Omega_\N: X_{m;N}(A) \to
X_m(A)\ {\rm as}\ N\to\infty\}\right) \ =  \ 1; \ee

\item (In probability) For each $m$, $X_{m;N} \to X_m$ in
probability if for all $\gep>0$, \be
\lim_{N\to\infty}\p_\N(|X_{m;N}(A) - X_m(A)|
> \gep) \ = \ 0;\ee

\item (Weak convergence) For each $m$, $X_{m;N} \to X_m$ weakly if
\be \p_\N(X_{m;N}(A) \le x) \ \to \ \p(X_m(A) \le x) \ee as
$N\to\infty$ for all $x$ at which $F_{X_m}(x) = \p(X_m(A) \le x)$ is
continuous.

\een

Alternate notations are to say \emph{with probability 1} for almost
sure convergence and \emph{in distribution} for weak convergence;
both almost sure convergence and convergence in probability imply
weak convergence. For our purposes we take $X_m$ as the random
variable which is identically $M_m$ (thus $X_m(A) = M_m$ for all $A
\in \Omega_\N$).

Our main tool to understand the $F^{A_N/\sqrt{N}}$ is the Moment
Convergence Theorem (see \cite{Ta} for example); our analysis is
greatly simplified by the fact that we have convergence to the
standard normal.

\begin{thm}[Moment Convergence Theorem]\label{thm:momct} Let $\{F_N(x)\}$ be a
sequence of distribution functions such that the moments \be M_{m;N}
\ = \ \int_{-\infty}^\infty x^m dF_N(x) \ee exist for all $m$. Let
$\Phi$ be the distribution function of the standard normal (whose
$m$\textsuperscript{{\rm th}} moment is $M_m$). If
$\lim_{N\to\infty} M_{m,N} = M_m$ then $\lim_{N\to\infty} F_N(x) =
\Phi(x)$. \end{thm}

\begin{defi}[Limiting spectral distribution]\label{defi:lsd}
If as $N\to\infty$ we have $F^{A_N/\sqrt{N}}$ converges in some
sense (for example, weakly or almost surely) to a distribution $F$,
then we say $F$ is the limiting spectral distribution of the
ensemble.
\end{defi}

In \S\ref{sec:intromainresults} we state our main results about the
type of convergence of the $F^{A_N/\sqrt{N}}$. The limiting spectral
distribution will be the distribution function of the standard
normal. The analysis proceeds by examining the convergence of the
moments. For example, assume for each $m$ that we have $X_{m,N}(A)
\to M_m$ almost surely. If \be B_m \ = \ \{A\in\Omega_\N: X_{m;N}(A)
\not\to M_m\ {\rm as}\ N\to\infty\},\ee then $\p(B_m) = 0$ and thus
\be \p\left(\bigcup_{m=0}^\infty B_m\right)\ =\ 0. \ee This and the
Moment Convergence Theorem allow us to conclude that with
probability 1, $F^{A_N/\sqrt{N}}(x)$ converges to $\Phi(x)$.

\subsection{Main Results}\label{sec:intromainresults}


By analyzing the moments of the $\mu_{A_N}$ (for $A_N$ an $N\times
N$ real symmetric palindromic Toeplitz matrix), we obtain results on
the convergence of $F^{A_N/\sqrt{N}}$ to the distribution function
of the standard normal. The $m$\textsuperscript{th} moment of
$\mu_{A_N}(x)$ is \be M_m(A_N) \ = \ \int_{-\infty}^\infty x^m
\mu_{A_N}(x)dx\ = \ \frac{1}{N^{\frac{m}{2}+1}} \sum_{i=1}^N
\gl_i^m(A_N). \ee

\begin{defi}\label{defi:MmN}
Let $M_m(N)$ be the average of $M_m(A_N)$ over the ensemble, with
each $A$ weighted by its distribution.  Set $M_m = \lim_{N\to\infty}
M_m(N)$. We call $M_m(N)$ the average $m^{\rm th}$ moment, and $M_m$
the limit of the average $m^{\rm th}$ moment.
\end{defi}

While we have two different definitions of $M_m$ (we have defined it
as both the $m$\textsuperscript{th} moment of the standard Gaussian
as well as the limit of $M_m(N)$), in Theorem \ref{thm:moments} we
prove that the $M_m(N)$ converge to the moments of the standard
Gaussian density, independent of $p$. Thus the two definitions are
the same. Specifically, $\lim_{N\to\infty} M_m(N) = (2k-1)!!$ if
$m=2k$ is even, and $0$ otherwise. Once we show this, then the same
techniques used in \cite{HM} allow us to conclude

\begin{thm}\label{thm:main} The limiting spectral
distribution of real symmetric palindromic Toeplitz matrices whose
independent entries are independently chosen from a probability
distribution $p$ with mean $0$, variance $1$ and finite higher
moments, converges weakly, in probability and almost surely to the
cumulative distribution function of the standard Gaussian, independent of $p$.
\end{thm}

We sketch the proof, which relies on Markov's method of moments.
While this technique has been replaced by other methods (which do
not have as stringent requirements on the underlying distribution),
the method of moments is well suited to random matrix theory
problems, as well as many questions in probabilistic number theory
(see \cite{Ell}).

By the eigenvalue trace lemma, \be\label{eq:MkANl2} \sumin \gl_i^m \
= \ \text{Trace}(A_N^m) \ = \ \sum_{1 \leq i_{1},\dots, i_{m} \leq
N}a_{i_{1} i_{2}}a_{i_{2} i_{3}}\cdots a_{i_{m} i_{1}}. \ee Applying
this to our palindromic Toeplitz matrices, we have \be\label{eqtr}
M_m(N) \ = \ \mathbb{E}[M_m(A_N)]\ = \ \frac{1}{N^{\frac{m}{2}+1}}
\sum_{1 \leq i_{1},\dots, i_m \leq
N}\E(b_{|i_{1}-i_{2}|}b_{|i_{2}-i_{3}|}\cdots b_{|i_{m}-i_{1}|}),
\ee where by $\E(\cdots)$ we mean averaging over the $N\times N$
palindromic Toeplitz ensemble with each matrix $A_N$ weighted by its
probability of occurring; thus the $b_j$ are i.i.d.r.v. drawn from
$p$. We show in \S\ref{sec:calcmoments} that the  $M_m =
\lim_{N\to\infty} M_m(N)$ are the moments of the standard Gaussian
density. The odd moment limits are easily shown to vanish, and the
additional symmetry (the palindromic condition) completely removes
the obstructions to the system of Diophantine equations studied in
\cite{HM}.

Convergence in probability follows from \be \lim_{N \to \infty}
\left(\E[ M_m(A_N)^2 ] - \E[M_m(A_N)]^2 \right) \ = \ 0, \ee
Chebyshev's inequality and the Moment Convergence Theorem, while
almost sure convergence follows from showing
\be\label{eq:thmprestrong} \lim_{N\to \infty} \E\left[ |M_m(A_N) -
\E[M_m(A_N)]|^4 \right] \ = \ O_m\left(\frac{1}{N^2}\right), \ee and
then applying Chebyshev's inequality, the Borel-Cantelli Lemma and
the Moment Convergence Theorem. Analogues of these estimates are
proven in \cite{HM} for the ensemble of real symmetric Toeplitz
matrices by degrees of freedom arguments concerning the tuples
$(i_1, \dots, i_m)$. The palindromic structure does not change the
number degrees of freedom, merely the contribution from each case.
Thus the arguments from \cite{HM} are applicable, and yield both
types of convergence. We sketch these arguments in
\S\ref{sec:weaksure}. In \S\ref{sec:conncircensemble} we investigate
related ensembles. In particular, we show our techniques apply to
real symmetric palindromic Hankel matrices, with Theorem
\ref{thm:main} holding for this ensemble as well. Further, we show
that the limiting spectral distribution of the palindromic Toeplitz
ensemble is the same as that of Bose and Mitra's symmetric Toeplitz
ensemble, implying that our Diophantine analysis is equivalent to
their analysis of the explicit formulas for the eigenvalues of their
ensemble. 

One particularly nice application of the correspondence between
these two ensembles is that we obtain a version of the almost sure
Central Limit Theorem for certain weighted sums of independent
random variables. Specifically, in \S\ref{sec:almostsureCLT} we show

\begin{thm}\label{thm:asclt2} Let $X_1, X_2, \dots$ be independent,
identically distributed random variables from a distribution $p$
with mean 0, variance 1, and finite higher moments. For $\omega =
(x_1,x_2,\dots)$ set $X_\ell(\omega) = x_\ell$, and consider the
probability space $(\Omega,\f,\p)$ (where $\p$ is induced from ${\rm
Prob}(X_\ell(\omega) \le x) = \int_{-\infty}^{x}p(t)dt$). Let \be
S_n^{(k)}(\omega) \ = \ \frac1{\sqrt{n/2}} \sum_{\ell=1}^{n}
X_\ell(\omega) \cos(\pi k \ell / n). \ee Then 
\be \p\left(\left\{\omega \in \Omega: \sup_{x\in\R}\left|\frac1n
\sum_{k=1}^{n} I_{S_n^{(k)}(\omega) \le x} - \Phi(x)\right| \to 0\
{\rm as}\ n \to \infty\right\}\right) \ = \ 1; \ee  here $I$ denotes
the indicator function and $\Phi$ is the distribution function of
the standard normal: \be \Phi(x) \ = \
\frac1{\sqrt{2\pi}}\int_{-\infty}^x e^{-t^2/2} dt.\ee
\end{thm}

 We
conclude in \S\ref{sec:futurework} by investigating the spacings
between normalized eigenvalues of palindromic Toeplitz matrices.


\section{Calculating the Moments}\label{sec:calcmoments}

Many of the calculations below are similar to ones in \cite{HM}, the
difference being that the additional symmetries imposed by the
palindromic condition remove the obstructions to the Diophantine
equations. Our main result, needed for the proof of Theorem
\ref{thm:main}, is that

\begin{thm}\label{thm:moments} For the ensemble of real symmetric
palindromic Toeplitz matrices with independent entries chosen
independently from a probability distribution $p$ with mean $0$,
variance $1$ and finite higher moments, each $M_m$ (the limit of the
average moments of the normalized empirical spectral measures)
equals the $m^{\rm th}$ moment of the standard Gaussian density.
Specifically, $M_{2k+1} = 0$ and $M_{2k} = (2k-1)!!$, where
$(2k-1)!! = (2k-1) \cdot (2k-3) \cdots 3 \cdot 1$.
\end{thm}

We prove Theorem \ref{thm:moments} in stages. In
\S\ref{sec:oddmoments} we show that the odd moments vanish, and that
the limit of the average zeroth and second moments are $1$.
Determining the moments is equivalent to counting the number of
solutions to a system of Diophantine equations. In
\S\ref{sec:highermoments} we prove some properties of the
Diophantine system of equations, which we then use in
\S\ref{sec:fourth} to show that $M_4$, the limit of the average
fourth moment as $N\to\infty$, equals that of the standard Gaussian
density. As we can always translate and rescale a probability
distribution with finite moments to have mean $0$ and variance $1$,
the first moment that shows the shape of an even distribution is the
fourth. This supports the claim that the palindromic condition
removes the obstructions. We then use linear algebra techniques (and
the ability to solve several Diophantine equations at once) to show
that the limits of all the even average moments agree with those of
the standard Gaussian density in \S\ref{sec:evenmoment}.

We introduce some notation. Let $A_N$ be an $N\times N$ real
symmetric palindromic Toeplitz matrix. We write $a_{ij}$ for the
entry in the $i$\textsuperscript{th} row and $j$\textsuperscript{th}
column. We determine which entries are forced to have the same value
as $a_{i_m i_{m+1}}$. As $A_N$ is a real symmetric palindromic
Toeplitz matrix, if $a_{i_n i_{n+1}}$ is forced to have the same
value then either (1) it is on the same diagonal; (2) it is on the
diagonal obtained by reflecting the diagonal $a_{i_m i_{m+1}}$ is on
about the main diagonal; (3) it is on the diagonal corresponding to
$b_{(N-1)-|i_{m+1}-i_m|}$; (4) it is on the diagonal obtained by
reflecting about the main diagonal the diagonal corresponding to
$b_{(N-1)-|i_{m+1}-i_m|}$. In other words, \be a_{i_m i_{m+1}} \ = \
a_{i_n i_{n+1}} \ {\rm if}\
\begin{cases} |i_{m+1}-i_m|\ =\ |i_{n+1}-i_n| \\ |i_{m+1}-i_m|\ =\ N-1 - |i_{n+1}-i_n|,
\end{cases} \ee where we set $i_{N+1}$ equal to $i_1$.
Equivalently, \be\label{eq:matchingdiags} a_{i_m i_{m+1}} \ = \
a_{i_n i_{n+1}} \ {\rm if}\
\begin{cases} i_{m+1} - i_m \ = \ \pm (i_{n+1} - i_n) \\ i_{m+1} - i_m \ = \ \pm
(i_{n+1} - i_n) + (N - 1) \\ i_{m+1} - i_m \ = \ \pm (i_{n+1} - i_n)
- (N - 1).
\end{cases} \ee

We denote the common value by $b_{|i-j|}$, and use $b_\alpha$ to
refer to a generic diagonal (thus the $a_{ij}$'s refer to individual
entries and the $b_\alpha$'s refer to diagonals). Each such matrix
is determined by choosing $\frac{N}2$ numbers independently from
$p$, the $b_\alpha$ with $\alpha \in \{0,\dots, \frac{N}2-1\}$. The
moments are determined by analyzing the expansion for $M_k(N)$ in
\eqref{eqtr}. We let $p_k$ denote the $k$\textsuperscript{th} moment
of $p$, which is finite by assumption.

We often use big-Oh notation: if $g(x)$ is a non-negative function
then $f(x) = O(g(x))$ (equivalently, $f(x) \ll g(x)$) if there are
constants $x_0, C > 0$ such that for all $x \ge x_0$, $|f(x)| \le
Cg(x)$. If the constant depends on a parameter $m$ we often write
$\ll_m$ or $O_m$.

\subsection{Zeroth, Second and Odd Moments}\label{sec:oddmoments}

\begin{lem}\label{oddthm}
Assume $p$ has mean zero, variance one and finite higher moments.
Then $M_0 = 1$ and $M_2 = 1$.
\end{lem}

\begin{proof}
For all $N$, $M_0(A_N) = M_0(N) = 1$. For the second moment, we have
\bea M_2(N) & \ = \ & \frac{1}{N^2} \sum_{1 \le i_1,i_2 \le N} \E(
a_{i_1 i_2}\cdot a_{i_2 i_1}) \nonumber\\ & \ = \ & \frac{1}{N^2}
\sum_{1 \le i_1,i_2 \le N} \E(a_{i_1 i_2}^2) \ = \ \frac{1}{N^2}
\sum_{1 \le i_1,i_2 \le N} \E(b_{|i_1-i_2|}^2). \eea As we have
drawn the $b_\alpha$'s from a variance $1$ distribution, the
expected value above is $1$.  Thus $M_2(N) = \frac{N^2}{N^2} = 1$,
so $M_2 = \lim_{N\to\infty} M_2(N)= 1$ also. \end{proof}

Note there are two degrees of freedom. We can choose $a_{i_1 i_2}$
to be on any diagonal. Once we have specified the diagonal, we can
then choose $i_1$ freely, which now determines $i_2$.

\begin{lem}\label{thm:odd}
Assume $p$ has mean zero, variance one and finite higher moments.
Then $M_{2k+1} = 0$. \end{lem}

\begin{proof} For $m=2k+1$ odd, in \eqref{eqtr} at least one $b_\alpha$
occurs to an odd power. If a $b_\alpha$ occurs to the first power,
as the expected value of a product of independent variables is the
product of the expected values, these terms contribute zero. Thus
the only contribution to an odd moment come when each $b_\alpha$ in
the expansion occurs at least twice, and at least one occurs three
times.

There are at most $k+1$ degrees of freedom. There are at most $k$
values of $b_\alpha$ to specify, and then once any index $i_\ell$ is
specified in \eqref{eqtr}, there are at most $8$ values (coming from
the four possible diagonals in \eqref{eq:matchingdiags}) for each
remaining index. Therefore of the $N^{2k+1}$ tuples $(i_1, \dots,
i_{2k+1})$, there are only $O(N^{k+1})$ tuples where the
corresponding $b_\alpha$'s are matched in at least pairs.

Consider such a tuple. Assume there are $r \le k$ different
$b_\alpha$, say $b_{\alpha_1}, \dots, b_{\alpha_r}$, with
$b_{\alpha_j}$ occurring $n_j \ge 2$ times (and further at least one
$n_j \ge 3$). Such an $(i_1, \dots, i_{2k+1})$ tuple contributes
$\prod_{j=1}^r \E\left[b_{\alpha_j}^{n_j}\right] = \prod_{j=1}^r
p_{n_j}$ to $M_{2k+1}(N)$, where $p_j$ is the
$j$\textsuperscript{th} moment of $p$ and hence finite. Thus this
term contributes $O_k(1)$ (where the constant depends on $k$); in
fact the constant is at most $\max_{j \le 2k+1} (|p_j|^k, 1)$.

Thus \bea M_{2k+1}(N)  \ \ll_k \ \frac{1}{N^{\frac{2k+1}{2}\ + \ 1}}
\cdot N^{k+1} \ \ll_k \ N^{-\foh}, \eea so $M_{2k+1} = \lim_{N \to
\infty} M_{2k+1}(N) = 0$, completing the proof.
\end{proof}

\subsection{Higher Moments}\label{sec:highermoments}

We expand on the method of proof of Lemma \ref{thm:odd} to determine
the even moments. We must find the $N\to\infty$  limit of
\be\label{eq:mknhmstart} M_{2k}(N) \ = \ \frac{1}{N^{k+1}}
\sum_{1\leq i_1,\dots,i_{2k} \leq N} \E (a_{i_1 i_2} a_{i_2 i_3}
\cdots a_{i_{2k} i_1}). \ee If the tuple $(i_1, \dots, i_{2k})$ has
$r$ different $b_\alpha$, say $b_{\alpha_1}, \dots, b_{\alpha_r}$,
with $b_{\alpha_j}$ occurring $n_j$ times, then the tuple
contributes $\prod_{j=1}^r \E\left[b_{\alpha_j}^{n_j}\right] =
\prod_{j=1}^r p_{n_j}$.

\begin{lem}\label{lem:mustbepairs} The tuples in
\eqref{eq:mknhmstart} where some $n_j \neq 2$ contribute
$O_{k}(\frac1{N})$ to $M_{2k}(N)$, the average
$2k$\textsuperscript{{\rm th}} moment. Thus, as $N\to\infty$, the
only tuples that contribute to $M_{2k}$ are those where the $a_{i_m
i_{m+1}}$ are matched in pairs.
\end{lem}

\begin{proof} If an $n_j = 1$ then the corresponding $b_{\alpha_j}$ occurs to
the first power. Its expected value is zero, and thus there is no
contribution from such tuples. Thus each $n_j \ge 2$, and the same
argument as in Lemma \ref{thm:odd} shows that each tuple's
contribution is $O_k(1)$. If an $n_j \ge 3$ then the corresponding
$b_{\alpha_j}$ occurs to the third or higher power, and there are
less than $k+1$ degrees of freedom (there are $O_k(N^k)$ tuples
where each $n_j \ge 2$ and at least one $n_j \ge 3$). As each
tuples' contribution is $O_k(1)$, and we divide by $N^{k+1}$ in
\eqref{eq:mknhmstart}, then the total contribution from these tuples
to $M_{2k}(N)$ will be $O_k(\frac1{N})$.  So in the limit as
$N\to\infty$ the contribution to $M_{2k}$ from tuples with at least
one $n_j \ge 3$ is $0$.
\end{proof}

\begin{rek}\label{freerek} Therefore the $b_{\alpha_j}$'s must be matched in
pairs.  There are $k+1$ degrees of freedom (we must specify values
of $b_{\alpha_1}, \dots, b_{\alpha_k}$, and then one index
$i_\ell$). It is often convenient to switch viewpoints from having
these $k$ pairings and one chosen index to having $k+1$ free indices
to choose, and we do so frequently. Another interpretation of Lemma
\ref{lem:mustbepairs} is that of the $N^{2k}$ tuples, only
$O_k(N^{k+1})$ have a chance of giving a non-zero contribution to
$M_{2k}(N)$. As any tuple contributes at most $O_k(1)$ to
$M_{2k}(N)$, in the arguments below we constantly use degree of
freedom arguments to show certain sets of tuples do not contribute
as $N\to\infty$ (specifically, any set of tuples of size $O_k(N^k)$
contributes $O_k(\frac1N)$ to $M_{2k}(N)$). \end{rek}

From \eqref{eq:matchingdiags}, if $a_{i_m i_{m+1}}$ is paired with
$a_{i_n i_{n+1}}$  then one of the following holds: \bea i_{m+1} -
i_m &\ =\ & \pm (i_{n+1} - i_n) \nonumber\\ i_{m+1} - i_m & = & \pm
(i_{n+1} - i_n) + (N - 1) \nonumber\\ i_{m+1} - i_m & = & \pm
(i_{n+1} - i_n) - (N - 1). \eea These equations can be written more
concisely. There is a choice of $C_\ell \in$ $\{0$, $\pm(N-1)\}$
($\ell$ is a function of the four indices) such that \be
\label{eqsgen} i_{m+1} - i_m \ = \ \pm (i_{n+1} - i_n) + C_\ell. \ee

The following lemma greatly prunes the number of possible matchings.

\begin{lem}\label{lem:allminussigns}
Consider all tuples $(i_1, \dots, i_{2k})$ such that the
corresponding $b_\alpha$'s are matched in pairs. The tuples with
some $a_{i_m i_{m+1}}$ paired with some $a_{i_n i_{n+1}}$ by a plus
sign in \eqref{eqsgen} contribute $O_k(\frac1N)$ to $M_{2k}(N)$.
Thus, as $N\to\infty$, they contribute $0$ to $M_{2k}$.
\end{lem}

\begin{proof}
Each tuple $(i_1, \dots, i_{2k})$ contributes \be \E[a_{i_1 i_2}
\cdots a_{i_{2k} i_1}] \ = \ \E[b_{|i_2-i_1|} \cdots b_{|i_1
-i_{2k}|}] \ee to $M_{2k}(N)$, and the only contributions we need
consider are when the $a_{i_m i_{m+1}}$ are matched in pairs. There
are $k$ equations of the form \eqref{eqsgen}; each equation has a
choice of sign (which we denote by $\gep_1, \dots, \gep_k$) and a
constant (which we denote by $C_1, \dots, C_k$; note each $C_\ell$
is restricted to being one of three values). We let $x_1, \dots,
x_k$ be the values of the $|i_{m+1} - i_{m}|$ on the left hand side
of these $k$ equations. Define $\widetilde{x}_1 = i_2 - i_1$,
$\widetilde{x}_2 = i_3 - i_2, \dots, \widetilde{x}_{2k} = i_{1} -
i_{2k}$.  We have
\bea\label{eq:relisxs} i_2 & \ = \ & i_1 - \widetilde{x}_1 \nonumber\\
i_3 & = & i_1 - \widetilde{x}_1 - \widetilde{x}_2 \nonumber\\ &
\vdots & \nonumber\\ i_1 & = & i_1 - \widetilde{x}_1 - \cdots -
\widetilde{x}_{2k}. \eea By the final relation for $i_1$, we find
\be\label{frelate} \widetilde{x}_1 + \cdots + \widetilde{x}_{2k} \ =
\ 0. \ee

Let us say that $a_{i_m i_{m+1}}$ is paired with $a_{i_n i_{n+1}}$.
Then we have relations between the indices because they must satisfy
one of the $k$ relations; let us assume they satisfy the
$\ell$\textsuperscript{th} equation. Further, by definition there is
an $\eta_\ell = \pm 1$ such that $i_{m+1} - i_{m} = \eta_\ell
x_\ell$; this is simply because we have defined the $x_j$'s to be
the absolute values of the $i_{m+1} - i_{m}$ on the left hand sides
of the $k$ equations. We therefore have that \be i_{m+1} - i_{m} \ =
\ \gep_\ell (i_{n+1} - i_{n}) + C_\ell, \ee or equivalently that \be
\widetilde{x}_m \ = \ \eta_\ell x_\ell \ = \ \gep_\ell
\widetilde{x}_n + C_\ell. \ee Since $\gep_\ell^2 = 1$, we have that
\be \widetilde{x}_n \ = \ \eta_\ell \gep_\ell x_\ell - \gep_\ell
C_\ell. \ee

Therefore each $x_\ell$ is associated to two $\widetilde{x}$'s, and
occurs exactly twice, once through $\widetilde{x}_m = \eta_\ell
x_\ell$ and once through $\widetilde{x}_n = \eta_\ell \gep_\ell
x_\ell - \gep_\ell C_\ell$. Substituting for the $\widetilde{x}$'s
in \eqref{frelate} yields \be\label{eq:csumtozero} \sum_{m=1}^{2k}
\widetilde{x}_m \ = \ \sum_{\ell=1}^{k} (\eta_\ell (1 + \gep_\ell )
x_\ell - \gep_\ell C_\ell )\ =\ 0. \ee

If any $\gep_\ell = 1$, then the $x_\ell$ are not linearly
independent, and we have fewer than $k+1$ degrees of freedom. There
will be at most $O_k(N^k)$ such tuples, each of which contributes at
most $O_k(1)$ to $M_{2k}(N)$. Thus the terms where at least one
$\gep_\ell = 1$ contribute $O_k\left(\fon\right)$ to $M_{2k}(N)$,
and are thus negligible in the limit. Therefore the only valid
assignment that can contribute as $N\to\infty$ is to have all
$\gep_\ell=-1$ (that is, only negative signs in \eqref{eqsgen}).
\end{proof}

\begin{rek}\label{rek:cssumzero}
The main term is when each $\gep_\ell = -1$. In this
case, \eqref{eq:csumtozero} immediately implies that the $C_\ell$'s
must sum to zero. This observation will be essential in analyzing
the even moments. \end{rek}

\subsection{The Fourth Moment}\label{sec:fourth}

We calculate the fourth moment in detail, as the calculation shows
how the palindromic structure removes the obstructions to the
Diophantine equations encountered in \cite{HM}.  This will establish
the techniques that we use to solve the general even moment in
\S\ref{sec:evenmoment}.

\begin{lem} Assume $p$ has mean zero, variance one and finite higher moments.
Then $M_4 = 3$, which is also the fourth moment of the standard
Gaussian density.
\end{lem}

\begin{proof} From \eqref{eqtr}, the proof follows by showing
\be\label{eq:M4} M_4 \ = \ \lim_{N\to\infty} \frac{1}{N^3}
\sum_{1\leq i,j,k,l \leq N} \E (a_{ij}a_{jk}a_{kl}a_{li}) \ee equals
$3$. From Lemma \ref{lem:mustbepairs}, the $a_{i_m i_{m+1}}$ must be
matched in pairs. There are three possibilities (see Figure
\ref{fig:4thconfigs}) for matching the $a_{i_m i_{m+1}}$ in pairs:
\bi
\item $(i,j)$ and $(j,k)$ satisfy \eqref{eqsgen}, and $(k,l)$ and
$(l,i)$ satisfy \eqref{eqsgen}; \item $(i,j)$ and $(k,l)$ satisfy
\eqref{eqsgen}, and $(j,k)$ and $(l,i)$ satisfy \eqref{eqsgen};
\item $(i,j)$ and $(l,i)$ satisfy \eqref{eqsgen}, and $(j,k)$ and
$(k,l)$ satisfy \eqref{eqsgen}. \ei

\begin{figure}
\begin{center}
\scalebox{.5}{\includegraphics{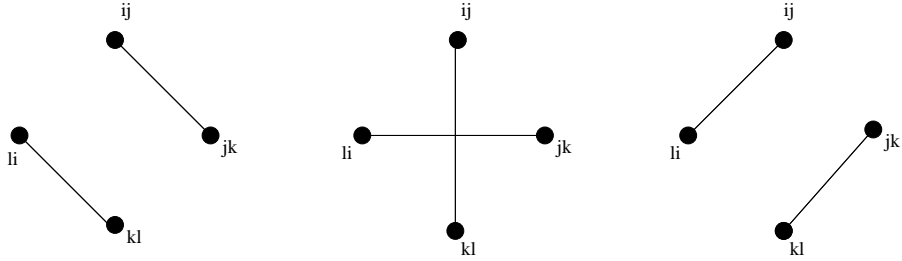}}
\caption{\label{fig:4thconfigs} Possible Configurations for the
Fourth Moment Matchings}
\end{center}\end{figure}

By symmetry (write $a_{ij}a_{jk}a_{kl}a_{li}$ as
$a_{li}a_{ij}a_{jk}a_{kl}$), the third case has the same
contribution as the first. These two cases are examples of adjacent
matchings. In the tuple $(i,j,k,l)$ we have four pairs, $(i,j),
(j,k), (k,l)$ and $(l,i)$, and we match the two adjacent ones. Note
that while in each case it is possible for both pairs to be
associated to the same $b_\alpha$ ($\alpha \in \{0, \dots,
\frac{N}2-1\}$), such tuples give a lower order contribution. We can
therefore ignore the contribution when both pairs have the same
value, as this is a correction of size $O(\frac1{N})$ to $M_4$.
Also, by Lemma \ref{lem:allminussigns}, we
only have minus signs in \eqref{eqsgen}.\\

\noindent \textbf{Case One: Adjacent Matching.} Consider the
adjacent matching (which occurs twice by relabeling). We thus have
the following pair of equations: \bea j - i \ = \ - (k - j ) + C_1,
\ \ \ \ \ \ \ \ \ \ l - k \ = \ - ( i - l ) + C_2 \eea Rewriting
these equations, we find that \bea k \ = \ i + C_1 \ \ \ \text{and}
\ \ \ k  \ =  \ i + C_2, \eea with $C_1, C_2 \in \{0, \pm(N-1)\}$
and $i,j,k,l \in \{1, \dots, N\}$.

We divide by $N^3$ in \eqref{eq:M4}. While we have $N^4$ tuples
$(i,j,k,l)$, only the $O(N^3)$ which have the $a_{i_m i_{m+1}}$
matched in pairs contribute. In fact, any set of tuples of size
$O(N^2)$ will not contribute in the limit. Thus we may assume $C_1$
and $C_2$ equal zero. For example, if $C_1 = N-1$ then $i$ is forced
to equal $1$, which forces $k$ to equal $N$. Letting $j$ and $l$
range over all possible values still gives only $N^2$ such tuples.
Similar arguments handle the case of $C_1 = -(N-1)$.

Thus $C_1 = C_2 = 0$; there are $N$ choices for $k \in
\{1,\dots,N\}$, and then $i$ is determined. We have $N$ choices for
$j \in \{1,\dots,N\}$ and $N-O(1)$ choices for $l$ (we want the two
pairs to correspond to different $b_\alpha$, so we must choose $l$
so that $a_{kl}$ is not on an equivalent diagonal to $a_{ij}$).
There are $N^3 - O(N^2)$ such tuples, each contributing $1$ (the
second moments of $p$ equal $1$, and we divide by $N^3$). Thus each
adjacent pairing case contributes $1 + O(\frac1N)$ to $M_4(N)$. As
there are two adjacent matching cases, as $N\to\infty$ these
contribute $2$ to $M_4$. \\

\noindent \textbf{Case Two: Non-adjacent Matchings.} The equations
for the non-adjacent case gives the following pair of equations:
\bea j - i \ = \ -(l-k) + C_1 \ \ \ \ \ \ \ \ \ \ k-j \ = \ -(i-l) +
C_2, \eea or equivalently \be j \ = \ i + k - l + C_1 \ = \ i + k -
l - C_2.\ee We see that $C_1 = -C_2$, or $C_1 + C_2 = 0$.

In \cite{HM}, as $N\to\infty$ this non-adjacent pairing contributed
$\frac23$ to $M_4$, and was responsible for the non-Gaussian
behavior. The difference is that in \cite{HM} we had the relation $j
= i+k-l$ \emph{without} the additional factor $C_1 \in \{0,
\pm(N-1)\}$. The problem was that we required each $i,j,k,l \in
\{1,\dots,N\}$; however, if we choose $i,k$ and $l$ freely then $j$
may not be in the required range. For example, whenever $i,k \ge
\frac23 N$ and $l < \frac13N$ then $j > N$; thus for the Toeplitz
ensemble at least $\frac1{27}N^3$ of the $N^3$ tuples that ``should
have'' contributed $1$ instead contributed $0$.

We now show this does not happen for the palindromic Toeplitz
ensemble. For any $i,k,l \in \{1,\dots, N\}$ there is a choice of
$C_1 \in \{0, \pm(N-1)\}$ such that $j\in \{1,\dots, N\}$ as well.
The choice of $C_1$ is unique unless $i+k-l \in \{1, N\}$, but this
is an additional restriction (i.e., we lose a degree of freedom
because an additional equation must be satisfied) and there are only
$O(N^2)$ triples $(i,k,l)$ with $i+k-l \in \{1,N\}$. Thus there are
again $N^3 + O(N^2)$ tuples, each with a contribution of $1$ (if all
four $a_{i_m i_{m+1}}$ are on equivalent diagonals then this is
again a lower order term, as there are at most $O(N^2)$ such
tuples). As there is one non-adjacent matching case, as $N\to\infty$
this
contribute $1$ to $M_4$. \\

Adding the contribution from the two cases gives a value of $3$ for
$M_4$, the limit of the average fourth moment, completing the proof.
\end{proof}

\subsection{The General Even Moment}\label{sec:evenmoment}

We now address the general case.  Using the linear algebra
techniques highlighted in the fourth moment calculation, we complete
the proof of Theorem \ref{thm:moments} by showing the limit of the
even average moments, the $M_{2k}$'s, agree with the even moments of
the standard Gaussian density.

Fix an even number $2k \ge 6$. By Lemma \ref{lem:mustbepairs} the
$a_{i_m i_{m+1}}$ must be matched in pairs. Each pair satisfies an
equation like \eqref{eqsgen}, and by Lemma \ref{lem:allminussigns}
the negative sign must hold. There are $(2k-1)!!$ ways to
match\footnote{There are $\ncr{2k}{2}$ ways to choose the first two
objects to be paired, $\ncr{2k-2}{2}$ ways to choose the second two
objects to be paired, and so on. As order does not matter, there are
$k!$ ways to arrange the $k$ pairs. Thus the number of matchings is
\be \ncr{2k}{2} \ncr{2k-2}{2} \cdots \ncr{2}{2}\ \Big/\ k! \ = \
(2k-1)!!.\nonumber\ \ee } the $2k$ objects in pairs. The proof of
Theorem \ref{thm:moments} is completed by showing that each of the
$(2k-1)!!$ matchings contributes $1$ to $M_{2k}$, as this then
implies that $M_{2k} = (2k-1)!!$.

Consider any matching of the $2k$ pairs of indices $(i_1, i_2 ),
(i_2, i_3), \dots, (i_{2k}, i_1)$ into $k$ pairs. We obtain a system
of $k$ equations. Each equation is of a similar form; for
definiteness we describe the equation when $(i_m, i_{m+1})$ is
paired with $(i_n, i_{n+1})$: \bea\label{gen1} i_{m+1} - i_m &\ =\ &
- (i_{n+1} - i_n) + C_j, \eea where as always $i_{2k+1} = i_1$, each
index is in $\{1,\dots,N\}$ and $C_j \in \{0, \pm(N-1)\}$. We may
re-write \eqref{gen1} as \bea\label{gen2} C_j &\ =\ & i_{m+1} - i_m
+ i_{n+1} -i_n. \eea Note that if we write the $k$ equations in the
form given by \eqref{gen2}, then each index $i_\alpha$ occurs
exactly twice. It occurs once with a coefficient of $+1$ and once
with a coefficient of $-1$. This is because the index $i_\ga$ occurs
in exactly two pairs of indices, in $(i_{\alpha-1}, i_{\alpha})$
(where it has a $+1$) and in $(i_\alpha, i_{\alpha+1})$ (where it
has a $-1$).

It is useful to switch between these two viewpoints (\eqref{gen1}
and \eqref{gen2}), and we do so below. We have $k$ equations and
$2k$ indices. We show there are $k+1$ degrees of freedom. In fact,
more is true. In the results that follow, we show $k+1$ of the
indices can be chosen freely in $\{1,\dots,N\}$, and for each
choice, there is a choice of the $C_j$'s such that there are values
for the remaining $k-1$ indices in $\{1,\dots,N\}$, and all $k$
equations hold. This means that each of the $N^{k+1}$ tuples (coming
from choosing $k+1$ of the indices freely) contributes $1$, which
shows this matching contributes $1$ to $M_{2k}$.

We first show how to determine which $k+1$ of the $2k$ indices we
should take as our free indices. Determining a good, general
procedure for finding the right free indices for an arbitrary choice
of the $(2k-1)!!$ matchings was the hardest step in the proof.

\begin{lem}\label{choicelem}
Consider the system of $k$ equations above, where each is of the
form described in \eqref{gen2}. We may number the equations from $1$
to $k$ and choose $k+1$ indices to be our free indices in such a way
that only the last equation has no dependent indices occurring for
the first time. For the first $k-1$ equations, there is always a
dependent index occurring for the first time, and there is always a
choice of the $C_j$'s so that the dependent indices in the first
$k-1$ equations take on values in $\{1,\dots, N\}$.
\end{lem}

It is important that in each equation only one dependent index
occurs for the first time. The reason is that we are trying to show
$N^{k+1} + O_k(N^k)$ of the $N^{k+1}$ choices of the independent
indices lead to valid configurations. If there were an equation with
dependent indices whose values were already determined, then we
would have restrictions on the independent indices and a loss of
degrees of freedom. We shall handle the last equation later (as
clearly every index occurring in the last equation has occurred in
an earlier equation).

\begin{proof} Choose any of the $k$ equations. We shall refer to it as
eq($k$).  This equation contains exactly four indices. As this is
the \emph{last} equation, each index must have appeared in an
earlier equation. Thus, eq($k$) marks the second time we have seen
each of these four indices.

Choose any of the four indices, and select the equation in which
this index first appeared. There is only one such equation, as each
index occurs in exactly two equations. We label this equation
eq($k-1$), and we let the index which we have just chosen be one of
our $k-1$ dependent indices. For the other three indices, either two
have a plus sign and the third has a negative sign (in which case
our dependent index has a negative sign), or two have a negative
sign and one has a positive sign (in which case our dependent index
has a positive sign). Let us assume our dependent index has a
negative sign, and consider the corresponding equation in the form
of \eqref{gen2}; the case where the dependent index has a positive
sign is handled similarly. The other three indices' sum is in
$\{2-N,\dots,2N-1\}$. If the sum is in $\{2-N,\dots,1\}$ we may take
$C_{k-1} = -(N-1)$; if the sum is in $\{1,\dots,N\}$ we may take
$C_{k-1} = 0$; if the sum is in $\{N,\dots,2N-1\}$ we may take
$C_{k-1} = N-1$. In each case there is a valid choice of the
dependent index. While if the other indices sum to $1$ or $N$ then
there are two choices of $C_{k-1}$, we shall see in Lemma
\ref{detlem} that this give lower order contributions and may be
safely ignored as $N\to\infty$.

Now consider the indices in eq($k$) and eq($k-1$); as long as at
least one index has appeared only once in these two equations, we
may continue the process. We choose any such index. It will be one
of our dependent indices, and we label the unique other equation it
occurs in as eq($k-2$).

We claim we may repeat this process until we have chosen one index
from all but eq($k$) as a dependent index, and each equation has a
dependent index which occurs for the first time in that equation.
The only potential problem is there is an $m
> 1$ such that, after we chose which equation to label eq($m$),
every index in eq($m$) through eq($k$) occurs exactly twice. If this
were so, we would not be able to continue and choose a new dependent
index and a new equation to be eq($m-1$). We show that there is no
such $m>1$.

We prove this by contradiction. Assume not, so every index in
eq($m$) through eq($k$) occurs twice. In our initial configuration,
we had $2k$ pairs of indices: $(i_1, i_2)$, $(i_2, i_3)$, $\dots$,
$(i_{2k}, i_1)$; note that each index is in exactly two pairs.
Without loss of generality, assume eq($k$) has index $i_1$. Our
assumptions imply we have both $i_1$'s, which means we have the
pairs $(i_1, i_2)$ and $(i_{2k}, i_1)$. Since we are assuming each
index which occurs, occurs twice, we have the other $i_2$ and the
other $i_{2k}$. Thus we have the pairs $(i_2, i_3)$ and $(i_{2k-1},
i_{2k})$. Continuing in this manner, for $m>1$ we see that if we
were to terminate at some equation eq($m$), then there would be at
least two indices occurring only once.

Therefore the process never breaks down. We may choose a labeling of
the remaining $k-1$ equations such that, in each equation, there is
one and only one new dependent index occurring for the first time.
The remaining $k+1$ indices are our free indices.
\end{proof}

We now have $k+1$ free indices, and $k-1$ dependent indices. There
are $N^{k+1}$ choices for the $k+1$ free indices. We show that,
except for $O_k(N^k)$ ``bad'' choices of indices, there are unique
choices for the dependent indices and the $C_j$'s such that all $k$
equations are satisfied, and all indices are in $\{1,\dots,N\}$. As
the contributions to $M_{2k}(N)$ are divided by $N^{k+1}$, the
$O_k(N^k)$ bad indices contribute $O_k(\frac1N)$ to $M_{2k}(N)$, and
the $N^{k+1} - O_k(N^k)$ ``good'' indices contribute $1 +
O_k(\frac1N)$. Thus the contribution to $M_{2k}$ from this matching
is $1$.

\begin{lem}\label{detlem}
Except for $O_k(N^k)$ choices of the $k+1$ free indices, all the
constants $C_j$ ($j \in \{1,\dots,k-1\}$) are determined uniquely in
the set $\{0, \pm(N-1)\}$, the dependent indices are uniquely
determined in $\{1,\dots,N\}$, and the first $k-1$ equations are
satisfied.
\end{lem}

\begin{proof}
By Lemma \ref{choicelem}, each of the first $k-1$ equations
determines a single dependent index. Consider the sum of the other
three indices in these equations. In proving Lemma \ref{choicelem}
we showed that the $C_j$'s are unique whenever these sums are not
$1$ or $N$, and whenever a $C_j$ was unique it lead to a unique
choice of the dependent index in $\{1,\dots,N\}$ such that the
equation was satisfied. If the sum were either of these values, this
would give us another equation, and a loss of at least one degree of
freedom. This is immediate if one of the three indices is an
independent index; if all are dependent indices, then we simply
substitute for them with independent indices, and obtain an equation
involving many indices, at least one of which is independent. Thus
we again gain a relation among our independent indices. There are
therefore $O_k(N^k)$ choices of the $k+1$ free indices such that the
$C_j$'s are not uniquely determined.
\end{proof}

Notice how in the previous lemma, the last coefficient $C_k$ is not
included.  This is because in the above lemma we absolutely needed
to be able to determine our dependent index (which occurred for the
first time in eq($j$)) with $C_j$. However, in the last equation,
all the indices are determined. We therefore cannot determine $C_k$
in quite the same way as we did for the other $C_j$'s. We now show
that there is a valid choice of $C_k$ for $N^{k+1} + O_k(N^k)$ of
the choices of the free indices.

\begin{thm}\label{bigthm}
For any of the $(2k-1)!!$ matchings of the $2k$ pairs of indices
into $k$ pairs, there are $k+1$ free indices and $k-1$ dependent
indices. For all but $O_k(N^k)$ choices of the free indices, every
$C_j$ is in $\{0, \pm(N-1)\}$, and is uniquely determined.
Furthermore, we have $\sum_{j=1}^k C_j = 0$. Thus each matching
contributes $1 + O_k(\frac1N)$ to $M_{2k}(N)$, or equivalently
contributes $1$ to $M_{2k}$. Thus $M_{2k} = (2k-1)!!$, the $2k^{\rm
th}$ moment of the standard Gaussian density.
\end{thm}

\begin{proof} We have proved much of Theorem \ref{bigthm} in Lemmas
\ref{choicelem} and \ref{detlem}. What we must show now is that, for
all but $O_k(N^k)$ ``bad'' choices of the free indices, the last
equation is consistent. By our earlier results, we know all
equations but possibly the last are satisfied, all dependent indices
are in $\{1,\dots,N\}$, and for all but the ``bad'' choices of
indices, the $C_1, \dots, C_{k-1}$ are uniquely determined and in
$\{0, \pm(N-1)\}$.

Consider now the last equation, eq($k$). From \eqref{gen2} and the
fact that all indices are in $\{1,\dots,N\}$, we see that there is a
choice of $C_k \in \R$ such that eq($k$) holds. We must show that
$C_k \in \{0, \pm(N-1)\}$.

We first note that $C_k \in [2-2N, 2N-2]$. This is because each
index is in $[1,N]$, and in \eqref{gen2} two indices occur with a
positive sign and two with a negative sign.

We see that $C_k$ is a multiple of $N-1$ by adding the $k$ equations
(eq($1$) through eq($k$)). Each index occurs twice, once with a
negative sign and once with a positive sign, and each $C_j$ occurs
once with a positive sign. Thus \be C_1 + \cdots + C_k \ = \ 0; \ee
see also Remark \ref{rek:cssumzero}. As $C_1, \dots, C_{k-1}$ $\in$
$\{0, \pm(N-1)\}$, we obtain that $C_k$ is a multiple of $N-1$. As
$C_k \in [-2(N-1), 2(N-1)]$, we see that $C_k$ $\in$ $\{0$,
$\pm(N-1)$, $\pm2(N-1)\}$. We now show that $C_k = \pm2(N-1)$ for at
most $O_k(N^k)$ choices of the free indices.

Consider the case when $C_k=2N-2$; the other case is handled
similarly. For this to be true, in eq($k$) the two indices with
positive signs must equal $N$ and the two indices with negative
signs must be $1$. If this happens, we impose relations on previous
equations. Thus, just as in Lemma \ref{detlem}, we lose a degree of
freedom, and there are only $O_k(N^k)$ choices of the free indices
such that $C_k=2N-2$.

Therefore, $C_k \in \{0, \pm(N-1)\}$ and is determined uniquely
(except for at most $O_k(N^k)$ choices), and all $k$ equations are
satisfied with indices in $\{1,\dots,N\}$.
\end{proof}

This completes our proof that the limit of the average even moments,
the $M_{2k}$'s, agree with the even moments of the standard Gaussian
density.


\section{Convergence in Probability and Almost Sure Convergence}\label{sec:weaksure}

Showing the limit of the average moments agree with the standard
Gaussian's moments is the first step in proving Theorem
\ref{thm:main}. To complete the proof, we must show convergence in
probability and almost sure convergence (both of which imply weak
convergence). Fortunately, the arguments in \cite{HM} are general
enough to be immediately applicable for convergence in probability;
a small amount of additional work is needed for almost sure
convergence. We use the notation of \S\ref{sec:subnotation} and
\S\ref{sec:intromainresults} and state the minor changes needed to
apply the results of \cite{HM} to finish the proof. \emph{Note: in
\cite{HM} it is assumed that each $b_0 = 0$; by Lemma
\ref{lem:cantakemaindiagzero} we may assume $b_0 = 0$ without
changing the limiting spectral distribution of the ensemble.}

\subsection{Convergence in Probability}

Let $A$ an infinite sequence of real numbers and let $A_N$ be the
associated $N\times N$ real symmetric palindromic Toeplitz matrix.
Let $X_{m;N}(A)$ be the random variable which equals the
$m$\textsuperscript{th} moment of the measure associated to $A_N$
and let $M_m$ be the $m$\textsuperscript{th} moment of the standard
Gaussian. Set $X_m(A) = M_m$. We have $X_{m;N} \to X_m$ in
probability if for all $\gep>0$  \be \lim_{N\to\infty}
\p_\N\left(\left\{A \in \Omega_\N: |X_{m;N}(A) - X_m(A)|
> \gep\right\}\right) \ = \ 0. \ee By Chebyshev's inequality we have \be
\p_\N\left(\left\{A \in \Omega_\N: |X_{m;N}(A) - M_m(N)| >
\gep\right\}\right) \ \le \ \frac{\E[M_m(A_N)^2 ] -
\E[M_m(A_N)]^2}{\gep^2}. \ee As $M_m(N) - M_m \to 0$ as
$N\to\infty$, it suffices to show for all $m$ that \be \lim_{N \to
\infty} \left(\E[ M_m(A_N)^2 ] - \E[M_m(A_N)]^2 \right) \ = \ 0 \ee
and then apply the Moment Convergence Theorem (Theorem
\ref{thm:momct}).


By \eqref{eqtr} we have \bea \E[M_m(A_N)^2] & \ = \ &
\frac{1}{N^{m+2}} \sum_{1 \le i_1,\dots,i_m \le N} \nonumber\\ & & \
\ \ \ \ \ \times \  \sum_{1 \le j_1,\dots,j_m \le N} \E[
b_{|i_1-i_2|} \cdots
b_{|i_m-i_1|} b_{|j_1-j_2|} \cdots b_{|j_m-j_1|} ] \nonumber\\
\E[M_m(A_N)]^2 & = & \frac{1}{N^{m+2}} \sum_{1 \le i_1,\dots,i_m \le
N} \E[ b_{|i_1-i_2|} \cdots b_{|i_m-i_1|}] \nonumber\\ & & \ \ \ \ \
\ \times \ \sum_{1 \le j_1,\dots,j_m \le N}\mathbb{E}[b_{|j_1-j_2|}
\cdots b_{|j_m-j_1|} ]. \eea

There are two possibilities: if the absolute values of the
differences from the $i$'s are not on equivalent diagonals with
those of the $j$'s, then these contribute equally to $\E[ M_m(A_N)^2
]$ and $\E[M_m(A_N)]^2$. We are left with estimating the difference
for the crossover cases, when the value of an $i_\alpha -
i_{\alpha+1} = \pm (j_\beta - j_{\beta + 1})$. The proof of the
analogous result for the real symmetric Toeplitz ensembles in
\cite{HM} is done entirely by counting degrees of freedom, and
showing that at least one degree of freedom is lost if there is a
crossover. Such arguments are immediately applicable here, and yield
the weak convergence. All that changes is our big-Oh constants; the
important point to remember is that each $C_j \in \{0, \pm(N-1)\}$,
which means there are at most $3^{2m}$ configurations where we apply
the arguments of \cite{HM}.

\subsection{Almost Sure Convergence}

Almost sure convergence follows from showing that for each
non-negative integer $m$ that \be X_{m;N}(A) \ \to \ X_m(A) \ = \
M_m \ \ {\rm almost\ surely,} \ee and then applying the Moment
Convergence Theorem (Theorem \ref{thm:momct}). The key step in
proving this is showing that \be\label{eq:thmprestrongB} \lim_{N\to
\infty} \E\left[ |M_m(A_N) - \E[M_m(A_N)]|^4 \right] \ = \
O_m\left(\frac{1}{N^2}\right). \ee The proof is completed by three
steps. By the triangle inequality, \be |M_m(A_N) - M_m|\ \ \le \ \
|M_m(A_N) - M_m(N)|\ +\ |M_m(N) - M_m|. \ee As the second term tends
to zero, it suffices to show the first tends to zero for almost all
$A$.

Chebychev's inequality states that, for any random variable $X$ with
mean zero and finite $\ell$\textsuperscript{th} moment, \be
\text{Prob}(|X| \ge \gep) \ \le \ \frac{\E[|X|^\ell]}{\gep^\ell}.
\ee Note $\E[ M_m(A_N) - M_m(N)] = 0$, and following \cite{HM} one
can show the fourth moment of $M_m(A_N) - M_m(N)$ is
$O_m\left(\frac1{N^2}\right)$; we will discuss this step in greater
detail below. Then Chebychev's inequality (with $\ell = 4$) yields
\be \p_\N(|X_{m;N}(A) - X_m(A)| \ge \gep) \ \le \ \frac{\E[|M_m(A_N)
- M_m(N)|^4]}{\gep^4} \ \le \ \frac{C_m}{N^2 \gep^4}. \ee






The proof of almost sure convergence is completed by applying the
Borel-Cantelli Lemma and proving \eqref{eq:thmprestrongB}; we sketch
the proof below.

We assume $p$ is even for convenience (though see Remark 6.17 of
\cite{HM}). A careful reading of the proofs in \S6 of \cite{HM} show
that analogues of most of the results hold in the palindromic case
as well, as most of the proofs are simple calculations based on the
number of degrees of freedom. The only theorems where some care is
required are Theorems 6.15 (see equation (50)) and 6.16 (see
equation (51)). In those two theorems, more than just degree of
freedom arguments are used; however, the same equations are true for
each of our configurations, and thus analogues of these results hold
in the palindromic case as well, completing the proof of almost sure
convergence.


\section{Connection to Circulant and Other
Ensembles}\label{sec:conncircensemble}

We show how our analysis of the Diophantine equations associated to
the ensemble of real symmetric palindromic Toeplitz matrices may be
used to study the ensembles related to circulant matrices
investigated by Bose and Mitra, as well as other ensembles (for
example, real symmetric palindromic Hankel matrices). We conclude by
showing the two methods combine nicely to yield an almost sure
Central Limit Theorem.

We first state two needed results.

\begin{lem}[Cauchy's interlacing
property]\label{lem:cauchyinterlace} Let $A_N$ be an $N \times N$
real symmetric matrix and $B$ be the $(N-1)\times (N-1)$ principal
sub-matrix of $A_N$. If $\lambda_N  \ge \cdots \ge \lambda_1$
(respectively, $\lambda_{N-1}' \ge \cdots \ge \lambda_{1}'$) are the
eigenvalues of $A_N$ (respectively, $B_N$), then \be \lambda_N\ \ge\
\lambda_{N-1}'\ \ge\ \lambda_{N-1} \ \ge\ \lambda_{N-2}' \ \ge \
\cdots \ \ge\ \lambda_{2} \ \ge \ \lambda_{1}' \ \ge \ \lambda_1.
\ee
\end{lem}

For a proof, see \cite{DH}. For us, the important consequence of the
Cauchy interlacing property is the Rank Inequality (Lemma 2.2 of
\cite{Bai}):

\begin{lem}[Rank Inequality]\label{lem:rankineqBai}
Let $A_N$ and $B_N$ be $N\times N$ Hermitian matrices. Then \be
\sup_{x \in \R} |F^{A_N}(x) - F^{B_N}(x)| \ \le \ \frac{{\rm
rank}(A_N-B_N)}N. \ee
\end{lem}

To prove the equivalence of our methods with the direct analysis of
explicit formulas for eigenvalues, all we need is a simple
consequence of the Rank Inequality:

\begin{lem}[Special Case of the Rank Inequality]\label{lem:rankineq}
Let $A_N$ be an $N\times N$ real symmetric matrix with principal
$(N-1)\times (N-1)$ sub-matrix $B_{N-1}$. Then \be \sup_{x \in \R}
|F^{A_N}(x) - F^{B_{N-1}}(x)| \ \le \ \frac4N. \ee
\end{lem}

\begin{proof} We may extend $B_{N-1}$ to be an $N\times N$ real symmetric
matrix $B_N$ by setting all entries of $B_{N}$ in either the
$N$\textsuperscript{th} row or the $N$\textsuperscript{th} column
(but not both) equal to zero, and the entry $b_{NN}$ to any number
we wish. We may now apply Lemma \ref{lem:rankineqBai} to $A_N$ and
$B_N$ (with ${\rm rank}(A_N-B_N) \le 2$), and then note that the
spectral measures of $B_{N-1}$ and $B_N$ are close.
\end{proof}

We frequently use the rank inequality to show that two $N\times N$
matrices with common $(N-1)\times (N-1)$ principal sub-matrix have
empirical spectral measures differing by negligible amounts (as
$N\to\infty$).

\subsection{Circulant Ensembles}\label{subsec:circens}

Bose and Mitra (see page 9 of \cite{BM}) study what they call
symmetric Toeplitz matrices; their matrices are of the form
\be\label{eq:defsymmtoep} S_N \ = \ N^{-1/2}
\left(\begin{array}{cccccc}
x_{0}  &  x_{1}  & x_{2}  & \cdots & x_{N-2} & x_{N-1} \\
x_{1} &  x_{0}  & x_{1}  & \cdots & x_{N-3} & x_{N-2} \\
\vdots & \vdots  & \vdots & \ddots & \vdots &\vdots \\
            x_{N-1} &  x_{N-2} & x_{N-3} & \cdots & x_1 & x_0 \\
  \end{array}\right), \ \ \ \ x_{N-j} \ = \
  x_{j}.\ee
They prove the limiting spectral distribution exists for this
ensemble, and is the standard Gaussian. Their proof starts with
explicit formulas for the eigenvalues of the matrices in terms of
the matrix entries. The rest of the argument is similar to their
analysis of the empirical eigenvalue distribution of circulant
matrices. The normality of the limiting spectral distribution (i.e.,
it being the standard Gaussian) follows from a detailed analysis of
the eigenvalues, and requires several explicit computations. In our
analysis, the normality is a consequence of each matching
contributing fully, and allows us to avoid having to compute
detailed properties of the eigenvalues of the ensemble.

\begin{thm}\label{thm:symmtoepispaltoep} The ensembles of real symmetric
palindromic Toeplitz matrices (see \eqref{eq:defrsptmat}) and
symmetric Toeplitz matrices (see \eqref{eq:defsymmtoep}) have the
same limiting spectral distribution when the independent entries are
chosen from a distribution with mean $0$, variance $1$ and finite
higher moments.
\end{thm}

\begin{proof} The ensemble of $2N\times 2N$ real symmetric palindromic
Toeplitz matrices (see \eqref{eq:defrsptmat}) is almost, but not
quite, the same as the $(2N-1)\times (2N-1)$ symmetric Toeplitz
matrices studied by Bose and Mitra. The difference between the two
is that the symmetric Toeplitz matrices are $(2N-1)\times (2N-1)$
principal sub-matrices of the $2N\times 2N$ palindromic Toeplitz
matrices. By the rank inequality, as $N\to\infty$ the normalized
limiting spectral distributions converge to a common value; similar
arguments relate $(2N-1)\times (2N-1)$ palindromic Toeplitz and
$2N\times 2N$ symmetric Toeplitz ensembles. Thus solving either
ensemble is equivalent to solving the other. \end{proof}

In the ensemble of real symmetric Toeplitz matrices investigated in
\cite{HM}, the authors assumed $b_0 = 0$, as all $b_0$ does is shift
each normalized eigenvalue by $b_0/\sqrt{N}$; this will not affect
the limiting spectral distribution. Though the palindromic Toeplitz
matrices have $b_0$'s off the main diagonal, the following lemma
shows that we may again take $b_0 = 0$ without affecting the
limiting spectral distribution.

\begin{lem}\label{lem:cantakemaindiagzero}
The limiting spectral distribution of the ensemble of real symmetric
palindromic Toeplitz matrices, with the $b_i$ i.i.d.r.v from a
probability distribution with mean $0$, variance $1$ and finite
higher moments, is unchanged if we additionally require $b_0$ to
equal zero. \end{lem}

\begin{proof} If $b_0$ only occurred on the main diagonal, then its
only effect would be to shift each normalized eigenvalue by
$b_0/\sqrt{N}$, which is negligible in the limit. The argument thus
reduces to showing that the two other occurrences of $b_0$ (in the
upper right and lower left corners of our palindromic Toeplitz
matrices) have negligible effect on the distribution of the
normalized eigenvalues.

The proof follows by multiple applications of the rank inequality
(Lemma \ref{lem:rankineq}). Given an $N\times N$ real symmetric
palindromic Toeplitz matrix $A_N$ as in \eqref{eq:defrsptmat}, let
$A_N'$ be the matrix with entries $a_{ij}' = a_{ij}$, except for
$a_{1N}=a_{N1}=0$, and let $B_{N-1}$ be the $(N-1)\times (N-1)$
principal sub-matrix common to both $A_N$ and $A_N'$. Thus the
normalized empirical spectral measures of $A_N$ and $A_N'$ are both
within $4/N$ of that of $B_{N-1}$, and therefore differ from each
other by at most $8/N$. Let now $A_N''$ be the same matrix as $A_N'$
except with the main diagonal entries replaced by $0$. The
normalized eigenvalues (recall we divide by $\sqrt{N}$) of $A_N'$
and $A_N''$ differ by $b_0/\sqrt{N}$. Therefore \be
F^{A_N'/\sqrt{N}}(x)\ =\
F^{A_N''/\sqrt{N}}\left(x-\frac{b_0}{\sqrt{N}}\right).\ee Thus the
normalized empirical spectral distributions for $A_N$, $A_N'$ and
$A_N''$ all differ by a negligible amount as $N\to\infty$, so the
respective limiting spectral distributions of these three ensembles
converge to the same distribution.
\end{proof}

\begin{rek} For $i \in \{1,2\}$, assume $f_i(N) = o(N)$.
Arguing as in the proof of Lemma \ref{lem:cantakemaindiagzero}, we
see our results immediately extend to the limit of ensembles of
$N\times N$ real symmetric matrices where the upper left and lower
right blocks may be of size $f_1(N)$ and $f_2(N)$, and the main
diagonal block (of size $(N-f_1(N)-f_2(N)) \times
(N-f_1(N)-f_2(N))$) is real symmetric, palindromic and Toeplitz. We
use Lemma \ref{lem:rankineqBai} to compare the spectral measure of
the $N\times N$ matrix with that of the related $N\times N$ matrix
where we have set all entries in the first $f_1(N)$ rows and
columns, and the last $f_2(N)$ rows and columns equal to zero.
\end{rek}

\subsection{Hankel Matrices}

Our results hold for a wider class of matrices. Recall a Hankel
matrix is of the form \be H_N \ = \ \left(\begin{array}{ccccccc}
b_{N-1}  &  b_{N-2}  & b_{N-3}  & \cdots & b_2 & b_1 & b_{0} \\
b_{N-2} &  b_{N-3}  & b_{N-4}  & \cdots & b_1 & b_0 & b_{-1} \\
b_{N-3} &  b_{N-4} & b_{N-5}  & \cdots & b_0 & b_{-1} & b_{-2} \\
\vdots & \vdots  & \vdots &  & \vdots & \vdots & \vdots \\
            b_{2} &  b_{1} & b_{0} & \cdots & b_{-N+5} & b_{-N+4} & b_{-N+3} \\
            b_{1} &  b_{0} & b_{-1} & \cdots & b_{-N+4} & b_{-N+3} & b_{-N+2} \\
            b_{0} &  b_{-1} & b_{-2} & \cdots & b_{-N+3} & b_{-N+2} & b_{-N+1} \\
  \end{array}\right).
\ee Let $J_N$ be the $N\times N$ matrix which is zero everywhere
except on the anti-main diagonal, where the entries are $1$. For
example, $J_2 = \mattwo{0}{1}{1}{0}$. Note $J_N^2 = I_N$ (where
$I_N$ is the $N\times N$ identity matrix), and $H_N J_N$ is a
Toeplitz matrix. If additionally $b_k = b_{-k}$ then $H_N J_N$ is a
real symmetric Toeplitz matrix, and finally if the first row of
$H_N$ is a palindrome then $H_N J_N$ is a real symmetric palindromic
Toeplitz matrix. We shall call such $H_N$ (where $b_k = b_{-k}$ and
the first row is a palindrome) real symmetric palindromic Hankel
matrices.

There is a one-to-one correspondence between real symmetric
palindromic Toeplitz and Hankel matrices. A simple calculation shows
that if $(H_N,T_N)$ is such a pair, then \be H_N J_N \ = \ J_N H_N \
= \ T_N. \ee In particular, this implies \be H_N^2 \ = \ H_N I_N H_N
\ = \ H_N J_N J_N H_N \ = \ T_N^2, \ee and hence by induction we
have that \be H_N^{2k} \ = \ T_N^{2k}. \ee To show the spectral
measures attached to eigenvalues of real symmetric palindromic
Toeplitz matrices converge to the standard Gaussian, all we needed
was \eqref{eq:MkANl2}. There we saw the calculation depends solely
on the trace of the even powers of our matrices. As there is a
one-to-one correspondence, Theorem \ref{thm:main} holds for real
symmetric palindromic Hankel matrices as well.


\section{An almost sure Central Limit
Theorem}\label{sec:almostsureCLT}

We discuss how our results, combined with those of Bose and Mitra,
yield a version of the almost sure Central Limit Theorem. See
\cite{BC} (and the numerous references therein) for more details as
well as several examples of such theorems. We are grateful to the
referee for pointing out this application of our results.

Bose and Mitra analyze the distribution of the eigenvalues of the
symmetric Toeplitz ensemble (see \S\ref{subsec:circens}) by using
the explicit formulas for the eigenvalues. For these $N\times N$
circulant matrices with entries in $\{x_0, \dots, x_{N-1}\}$ (with
$x_{N-j} = x_j$), if $N$ is odd then the eigenvalues are \bea
\lambda_k \ = \ \frac1{\sqrt{N}} \sum_{\ell=0}^{N-1} x_\ell
\cos(2\pi k \ell / N) \ = \ \lambda_{N-k}; \eea a similar formula
holds if $N$ is even (there the eigenvalue $\lambda_{N/2}$ will have
multiplicity one). We have shown that the limiting distribution of
eigenvalues of this symmetric Toeplitz ensemble is the same as that
of our palindromic Toeplitz ensembles. The importance of this
connection is that we have shown the convergence is almost sure for
the palindromic ensemble. Thus we may translate this almost sure
convergence to a statement about the eigenvalues $\lambda_k$, which
are weighted sums of the symmetric Toeplitz matrix entries. We thus
obtain

\begin{thm}\label{thm:asclt1}
For each $N$ let $X_0, \dots, X_{N-1}$ be independent, identically
distributed random variables (subject to the condition that $X_j =
X_{N-j}$) from a distribution $p$ with mean 0, variance 1, and
finite higher moments. For $\omega = (x_0,x_1,\dots)$ set
$X_\ell(\omega) = x_\ell$, and consider the probability space
$(\Omega,\f,\p)$ (where $\p$ is induced from ${\rm
Prob}(X_\ell(\omega) \le x) = \int_{-\infty}^{x}p(t)dt$). Set \be
S_N^{(k)}(\omega) \ = \ \frac1{\sqrt{N}} \sum_{\ell=0}^{N-1}
X_\ell(\omega) \cos(2\pi k \ell / N). \ee Then 
\be \p\left(\left\{\omega \in \Omega: \sup_{x\in\R}\left|\frac1N
\sum_{k=0}^{N-1} I_{S_N^{(k)}(\omega) \le x} - \Phi(x)\right| \to 0\
{\rm as}\ n \to \infty\right\}\right) \ = \ 1; \ee here $I$ denotes
the indicator function and $\Phi$ is the distribution function of
the standard normal: \be \Phi(x) \ = \
\frac1{\sqrt{2\pi}}\int_{-\infty}^x e^{-t^2/2} dt \ee.
\end{thm}

A more useful version of the above is to note that we have double
counted all the eigenvalues (except possibly one, which will not
affect anything in the limit). Letting $N = 2n$ and looking at only
\emph{half} the eigenvalues, we immediately obtain\footnote{The
proof uses the fact that $S_N^{(k)}(\omega) = S_n^{(k)}(\omega)$,
which follows from our normalizations, simple algebra, and the fact
that $X_j = X_{N-j}$.} Theorem \ref{thm:asclt2} from Theorem
\ref{thm:asclt1}.


\section{Future Work}\label{sec:futurework}

So far we have investigated the density of the eigenvalues; we now
consider another problem, that of the spacings between adjacent
eigenvalues. Note the palindromic condition means that $0$ is always
an eigenvalue (because the first and last rows are identical),
though as $N\to\infty$, the contribution of one eigenvalue becomes
negligible.

As there are only $(N-2)/2$ degrees of freedom for the ensemble of
real symmetric palindromic Toeplitz matrices, which is much smaller
than $N(N+1)/2$, it is reasonable to believe the spacings between
adjacent normalized eigenvalues $\left(\frac{\gl_{i+1}(A)}{\sqrt{N}}
- \frac{\gl_i(A)}{\sqrt{N}}\right)$ may differ from those of full
real symmetric matrices. The ensemble of all real symmetric matrices
is conjectured to have normalized spacings given by the GOE
distribution (which is well approximated by $Axe^{-Bx^2}$) whenever
the independent matrix elements are independently chosen from a nice
distribution $p$. Studying thin sub-ensembles opens up the
possibility of seeing different behavior.

Interestingly (see \cite{JMRR} among others), the spacings between
adjacent normalized eigenvalues of $d$-regular graphs appear to be
given by the GOE as well. Thus, while the density of eigenvalues of
$d$-regular graphs is different than those of all real symmetric
matrices (Kesten's measure versus the semi-circle), the adjacent
normalized differences between eigenvalues behave like differences
of full real symmetric matrices. In the opposite extreme, consider
band matrices of width 1 (i.e., diagonal matrices). There the
spacing between adjacent normalized eigenvalues is Poissonian
($e^{-x}$), and the density of normalized eigenvalues is whatever
distribution the entries are drawn form.

\begin{figure}[h]
\begin{center}
\scalebox{.65}{\includegraphics{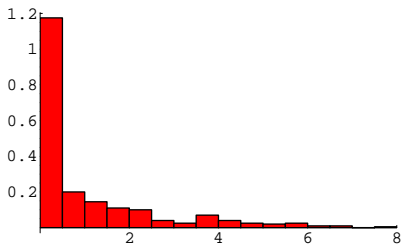}}
\caption{\label{fig:paltoepspac} Differences between normalized
eigenvalues 506 through 516 of 40 real symmetric palindromic $1000
\times 1000$ Toeplitz matrices.}
\end{center}
\end{figure}

\begin{figure}[h]
\begin{center}
\scalebox{.65}{\includegraphics{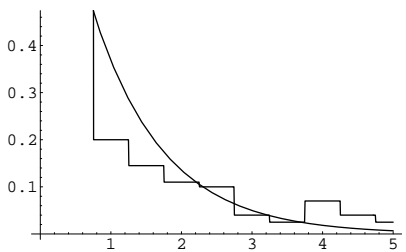}}
\caption{\label{fig:paltoepspacB} Comparison of differences between
normalized eigenvalues 506 through 516 of 40 real symmetric
palindromic $1000 \times 1000$ Toeplitz matrices and the standard
exponential; the small differences have been suppressed.}
\end{center}\end{figure}

We chose 40 Toeplitz matrices ($1000 \times 1000$) with entries
i.i.d.r.v. from the standard normal. The palindromic condition
implies that $0$ is always an eigenvalue of a real symmetric
palindromic Toeplitz matrix. To minimize the effect of this forced
eigenvalue, instead of looking at the middle 11 normalized
eigenvalues of each matrix, we looked at the next set of $11$
eigenvalues. This gave us $10$ differences between adjacent
normalized eigenvalues, and we compared those to the standard
exponential; if the spacings are Poissonian, the standard
exponential should be a good fit. Similar results were obtained for
larger shifts. See Figures \ref{fig:paltoepspac} and
\ref{fig:paltoepspacB} for the plots.

The distribution of differences looks approximately Poissonian;
definitely more Poissonian than GOE or GUE (both of which have small
probabilities of small spacings). While the fit to Poissonian
behavior is not as good as the real symmetric Toeplitz matrices
investigated in \cite{HM}, it is not unreasonable to conjecture that
in the limit as $N \to \infty$, the local spacings between adjacent
normalized eigenvalues will be Poissonian.


\section*{Acknowledgements} This work was performed at summer research programs
at The Ohio State University in 2004 and Brown University in 2005;
it is a pleasure to thank both institutions for their help and
support, as well as other program participants, especially Chris
Hammond, John Ramey and Jason Teich for many enlightening
discussions. We would also like to thank the anonymous referee for
very helpful comments on several drafts of this paper.


\ \\

\end{document}